\documentclass[4pt]{article} 

\usepackage[utf8]{inputenc} 
\usepackage{geometry} 
\usepackage{graphicx} 

\usepackage[colorlinks,citecolor=red]{hyperref}
\usepackage{amsmath}
\usepackage{amsfonts}
\usepackage{dsfont}
\usepackage{mathrsfs}
\usepackage{latexsym}
\usepackage{graphicx}
\usepackage{amssymb}
\usepackage{float}
\usepackage{epsfig}
\usepackage{epstopdf}
\usepackage{color,xcolor}
\usepackage{booktabs} 
\usepackage{array} 
\usepackage{paralist} 
\usepackage{verbatim} 
\usepackage{subfig} 

\newtheorem{Theorem}{Theorem}[section]
\newtheorem{Case}[Theorem]{Case}

\newtheorem{Corollary}[Theorem]{Corollary}

\newtheorem{Example}[Theorem]{Example}
\newtheorem{Lemma}[Theorem]{Lemma}

\newtheorem{Remark}[Theorem]{Remark}

\usepackage{fancyhdr} 
\usepackage{sectsty}
\allsectionsfont{\sffamily\mdseries\upshape} 

\usepackage[nottoc,notlof,notlot]{tocbibind} 
\usepackage[titles,subfigure]{tocloft} 

\def\0{\emptyset}

\def\p{{\bf Proof.~~}}
\def\q{\hfill\rule{1ex}{1ex}}

\title{Fractional matching number and spectral radius of \\ nonnegative matrix of graphs}

\author{Ruifang Liu \thanks{School of Mathematics and Statistics, Zhengzhou
University, Zhengzhou, Henan 450001, China. Email: rfliu@zzu.edu.cn},
Hong-Jian Lai \thanks{Corresponding author. Department of Mathematics, West Virginia
University, Morgantown, WV 26506, USA. E-mail: hjlai@math.wvu.edu},
Litao Guo \thanks{School of Applied Mathematics, Xiamen University
of Technology, Xiamen, Fujian 361024, China. Email: ltguo2012@126.com},
Jie Xue \thanks{School of Mathematics and Statistics, Zhengzhou
University, Zhengzhou, Henan 450001, China. Email: jie\_xue@126.com}}

\date{}

\begin{document}
\maketitle

\begin{abstract}
A fractional matching of a graph $G$ is a function $f:E(G) \to [0,1]$ such that for any $v\in V(G)$, $\sum_{e\in E_G(v)}f(e)\leq 1$
where $E_G(v) = \{e \in E(G): e$ is incident with $v$ in $G\}$. The fractional matching number of $G$
is $\mu_{f}(G) = \max\{\sum_{e\in E(G)} f(e): f$ is fractional matching of $G\}$.
For any real numbers $a \ge 0$ and $k \in (0, n)$, it is observed that
if $n = |V(G)|$ and $\delta(G) > \frac{n-k}{2}$, then $\mu_{f}(G)>\frac{n-k}{2}$.
We determine a function $\varphi(a, n,\delta, k)$ and show that
for a connected graph $G$ with $n = |V(G)|$, $\delta(G) \leq\frac{n-k}{2}$,
spectral radius $\lambda_1(G)$ and complement $\overline{G}$, each of the following holds.
\\
(i) If $\lambda_{1}(aD(G)+A(G))<\varphi(a, n, \delta, k),$ then $\mu_{f}(G)>\frac{n-k}{2}.$
\\
(ii) If $\lambda_{1}(aD(\overline{G})+A(\overline{G}))<(a+1)(\delta+k-1),$ then $\mu_{f}(G)>\frac{n-k}{2}.$
\\
As corollaries, sufficient spectral condition for fractional perfect matchings and
analogous results involving $Q$-index and $A_{\alpha}$-spectral radius
are obtained, and former spectral results in [European J. Combin. 55 (2016) 144-148] are extended.


\end{abstract}

\section{Introduction}
Graphs considered are simple and undirected. We generally follow \cite{BoMu08} for undefined terms and notation.
For a graph $G$, and a vertex $v \in V(G)$, define
$N_G(v) = \{u \in V(G): uv \in E(G)\}$ and $E_G(v) = \{e \in E(G): e$ is incident with $v$ in $G\}$.
If $S \subseteq V(G)$, then define $N_G(S) = \cup_{v \in S} N_G(v)$.
As in \cite{BoMu08}, we use $K_n$, $\delta(G)$, $\Delta(G)$, $\kappa(G)$ and $\overline{G}$ to denote the complete graph
of order $n$, the minimum degree, the maximum degree,
the connectivity and the complement of $G$, respectively.
For a subset $S\subseteq V(G)$ or $S \subseteq E(G)$,  $G[S]$ denotes the subgraph of $G$ induced by $S$.
A cycle (or path, respectively) on $n$ vertices is often denoted by $C_n$ (or $P_n$,
respectively). A cycle containing all vertices of a graph $G$ is called a Hamilton cycle of $G$.
Let $G$ be a
graph with vertex set $V(G)=\{v_{1},v_{2},\ldots,v_{n}\}$ and edge set $E(G)$. The {\bf adjacency matrix} of $G$
is the $n$ by $n$ matrix $A(G): = (a_{ij})$, where $a_{ij}$ is the number of edges joining $v_i$ and $v_j$ in $G$,
and the {\bf diagonal degree matrix} of $G$ is  the $n$ by $n$ matrix $D(G): = (d_{ij})$, where $d_{ij} = 0$ if $i \neq j$ and
$d_{ii} = d_G(v_i)$, for all $i$ with $1 \le i \le n$. For a subset $X \subseteq E(G)$, the
{\bf characteristic function} of $X$ is a function $f_X: E \to \{0, 1\}$ such that $f(e) = 1$ if and only if $e \in X$.

Our approach is motivated by the methods deployed in \cite{HGLL16, LiHL14, LHGL14, LiLT19, LiLT19a},
among others. For any graph $G$ with
the adjacency matrix $A(G)$ and the diagonal degree matrix $D(G)$ and for real numbers $a \ge 0$ and $b > 0$,
define $J(G: a,b) = aD(G)+bA(G)$. As $J(G:a,b)$ is also real and symmetric, all its eigenvalues are real.
Let $\lambda_i(J(G:a,b) )$ to be the $i$th largest eigenvalue of $J(G:a,b) $.
In particular, $\lambda_1(J(G:0,1)) = \lambda_1(G)$ is known as the {\bf spectral radius} of $G$;
$\lambda_1(J(G:1,1)) = q_1(G)$ is the {\bf $Q$-index} of $G$; and $\lambda_1(J(G: \alpha, 1-\alpha))=\lambda_1(A_{\alpha}(G))$ is called the
 {\bf $A_{\alpha}$-spectral radius} of $G$ with $\alpha\in[0,1]$, introduced by Nikiforov in \cite{Niki17}.

As in \cite{BoMu08}, an edge subset $M \subseteq E(G)$ is a {\bf matching} of a graph
$G$ if $\Delta(G[M]) \le 1$. The {\bf matching number} of $G$ is $\alpha'(G) = \max\{|M|: $
$M$ is a matching of $G\}$. A {\bf fractional matching} of a graph $G$ is a function
$f:E(G) \to [0,1]$ such that for any $v\in V(G)$, $\sum_{e\in E_G(v)}f(e)\leq 1$.
Thus if $f(e)\in\{0,1\}$ for every edge $e \in E(G)$, then $f$ is the characteristic function
of a matching in $G$. The {\bf fractional matching number} of $G$
is $\mu_{f}(G) = \max\{\sum_{e\in E(G)} f(e) : f$ is fractional matching of $G\}$.
The following is known.

\begin{Theorem}\label{2mu}
Let $G$ be a graph with $n = |V(G)|$. Each of the following holds.
\\
(i) (Lemma 2.1.2 of  \cite{ScUl97})  Any
fractional matching $f$ satisfies $\mu_{f}(G)\leq n/2$.
\\
(ii) (Theorem 2.1.3 of  \cite{ScUl97})
If $G$ is bipartite, then $\mu_f(G) = \alpha'(G)$.
\\
(iii) (Theorem 2.1.5 of  \cite{ScUl97})
$2 \mu_f(G)$ is an integer.
\end{Theorem}
A fractional matching $f$ of $G$
is a {\bf fractional perfect matching} if $\sum_{e\in E(G)} f(e)=n/2$. If a fractional perfect
matching takes values only in $\{0, 1\}$, then it is the characteristic function of a perfect
matching.

Let $i(G)$ be the number of isolated vertices of a graph $G$.
In \cite{ScUl97}, the following are displayed.
\begin{Theorem} \label{old1} (Scheinerman and Ullman \cite{ScUl97})
Let $G$ be a graph. Each of the following holds.
\\
(i) $\mu_{f}(G)\geq \alpha'(G)$.
\\
(ii) (The fractional Berge-Tutte Formula)
\begin{equation} \label{BT}
\mu_{f}(G)=\frac{1}{2}\left(n-\max\left\{i(G-S)-|S|: \forall S\subseteq V(G) \right\}\right).
\end{equation}
(iii) (The fractional Tutte's 1-Factor Theorem) A graph $G$ has a fractional perfect matching if and only if
\begin{equation} \label{T1}
\forall S\subseteq V(G), \; i(G-S)\leq|S|.
\end{equation}
\end{Theorem}
The investigation on the relationship between the eigenvalues and the matching number of a graph
was initiated by Brouwer and Haemers \cite{BrHa05}, in which sufficient conditions on $\lambda_{3}(G)$
to assure the existence of a perfect matching were discussed.
Cioab\v{a} et al. in \cite{Cioa05, CiGr07, CiGH09} improved and generalized these results
to obtain a best possible upper bound on $\lambda_{3}(G)$ for the existence of a perfect matching.
Furthermore, O and Cioab\v{a} \cite{OsCi10} determined the relationship between the eigenvalues of
a $t$-edge-connected $k$-regular graph and its matching number when $t\leq k-2$.

On the fractional matching number, O \cite{SO16} studied the connections between the fractional
matching number and the spectral radius of a connected graph
with given minimum degree. Xue et al. \cite{XuZS18} considered the relationship between
the fractional matching number and the Laplacian spectral radius of a graph.
The main goal of this study is to investigate the relationship between the fractional matching number
and the largest eigenvalue of the nonnegative matrix $J(G:a,b)$ of a graph $G.$
This provides a mechanism to have a unified approach to both adjacency eigenvalues and signless Laplacian eigenvalues,
and is quite different from those using only the Laplacian matrix. Consequently,
analogous results involving the $Q$-index and $A_{\alpha}$-spectral radius of a graph $G$ are also obtained, and
former results on spectral radius in \cite{SO16} are extended.

We first observe that (see Theorem \ref{th2} in the next section) for any real number $k\in(0, n)$
and a connected graph $G$ with $n=|V(G)|$ and  $\delta = \delta(G)$, if
$\delta>\frac{n-k}{2}$, then $\mu_{f}(G)>\frac{n-k}{2}$. Thus we focus on the investigation
of $\mu_{f}(G)$ when $\delta\leq\frac{n-k}{2}$.

For real numbers $a$, $k$ and an integer $\delta$, with $a \ge 0$,
$k\in(0, n)$ and $1\leq\delta\leq\frac{n-k}{2}$,  define
\begin{equation}\label{f}
\varphi(a, n,\delta, k)=
\left\{
\begin{array}{lc}
\delta\sqrt{1+\frac{2k}{n-k}} &\,~~~ \text{\mbox{if}~ $a=0$},
\\
\frac{2a\delta n}{n-k} &\,~~~ \text{\mbox{if}~ $a\in(0, 1]$},
\\
\frac{a\delta(n+k)}{n-k} &\,~~~ \text{\mbox{if}~  $a\in(1, +\infty)$}.
\end{array}
\right.
\end{equation}
Our main results are the following.

\begin{Theorem}\label{th1}
Let $a$ be a real number with $a\geq0$, $k\in(0, n)$ be a real number, and
$G$ be a connected simple graph with $n= |V(G)|$  and $\delta = \delta(G) \leq\frac{n-k}{2}.$ If
$$\lambda_{1}(aD(G)+A(G))<\varphi(a, n, \delta, k),$$ then $\mu_{f}(G)>\frac{n-k}{2}.$
\end{Theorem}

The upper bound in Theorem \ref{th1} is best possible when $a \in \{0,1\}$, as will be shown in Section 3.
A sufficient condition to ensure $\mu_{f}(G)>\frac{n-k}{2}$ in terms of $\lambda_{1}(J(\overline{G}: a, 1))$ is obtained
in Theorem \ref{th5} below. The upper bound in Theorem \ref{th5} is also best possible in some sense.

\begin{Theorem}\label{th5}
Let $a$ be a real number with $a\geq0$, $k\in(0, n)$ be a real number, and $G$ be a connected simple
graph with $n= |V(G)|$  and $\delta = \delta(G) \leq\frac{n-k}{2}$. If
$$\lambda_{1}(aD(\overline{G})+A(\overline{G}))<(a+1)(\delta+k-1),$$
then $\mu_{f}(G)>\frac{n-k}{2}.$
\end{Theorem}

The main results can be applied to obtain sufficient conditions in terms of $Q$-index
$q_1(G)$ and $A_{\alpha}$-spectral radius $\lambda_1(A_{\alpha}(G))$ for a simple graph $G$ to have a large value of $\mu_{f}(G)$.
Detailed justifications of the following corollaries will be presented in Section 4.

\begin{Corollary}\label{co3}
Let $G$ be a connected simple graph with $n = |V(G)|$ and $\delta = \delta(G) \leq\frac{n-k}{2}$. Each of the following holds.
\\
(i) If
$\displaystyle q_{1}(G)<\frac{2\delta n}{n-k}$, then $\mu_{f}(G)>\frac{n-k}{2}$.
\\
(ii) If $\displaystyle q_{1}(\overline{G})<2(\delta+k-1)$, then $\mu_{f}(G)>\frac{n-k}{2}$.
\end{Corollary}

\begin{Corollary}\label{co4}
Let $G$ be a connected simple graph with $n = |V(G)|$ and $\delta = \delta(G) \leq\frac{n-k}{2}$, and let $\alpha$
be a real number with $0 \le \alpha \le 1$. Each of the following holds.
\\
(i) If $\alpha=0$ and
$\displaystyle \lambda_1(A_{\alpha}(G))<\delta\sqrt{1+\frac{2k}{n-k}}$, then $\mu_{f}(G)>\frac{n-k}{2}$.
\\
(ii) If $0<\alpha\leq\frac{1}{2}$ and $\displaystyle \lambda_1(A_{\alpha}(G))<\frac{2\alpha\delta n}{n-k}$, then $\mu_{f}(G)>\frac{n-k}{2}$.
\\
(iii) If $\frac{1}{2}<\alpha\leq 1$ and $\displaystyle \lambda_1(A_{\alpha}(G))<\frac{\alpha\delta(n+k)}{n-k}$, then $\mu_{f}(G)>\frac{n-k}{2}$.
\\
(iv) If $\displaystyle \lambda_1(A_{\alpha}(\overline{G}))<\delta+k-1$, then $\mu_{f}(G)>\frac{n-k}{2}$.
\end{Corollary}

Preliminaries will be presented in the next section. Proofs of the main results and discussions of the
consequences will be given in the last two sections.

\section{Preliminaries}

Throughout this section, we always assume that $k\in(0, n)$
is a real number, and $G$ is a connected simple graph with $n= |V(G)|$ and $\delta=\delta(G)$.
A path with end vertices $u$ and $v$ is called a $(u,v)$-path.
We start with some of the well known results of Dirac.

\begin{Theorem} (see also Bondy and Murty \cite{BoMu08})\label{le1}
Let $G$ be a simple graph with $n = |V(G)|$ and $\delta=\delta(G) > 0$. Each of the following holds.
\\
(i) (Dirac \cite{Dira52}) If $\delta\geq n/2$ and $n\geq 3$, then $G$ is hamiltonian.
\\
(ii) (Dirac \cite{Dira52}) If $\kappa(G) \ge 2$ and $\delta\leq n/2$, then $G$ contains a cycle of length at least $2\delta$.
\\
(iii) $G$ contains a path of length $\delta$.
\end{Theorem}

Theorem \ref{le1} can be applied to prove Theorem \ref{th2} below, which indicates that the focus of the research should be restricted to
the cases when $\delta \le \frac{n-k}{2}$.

\begin{Theorem}\label{th2}
Let $s \ge 1$ be an integer. Each of the following holds.
\\
(i) If $G$ has a subgraph $L$ with $\Delta(L) \le 2$ and $|E(L)| \ge 2s$, then $\mu_{f}(G) \ge s$.
In particular, if $G$ has a cycle or path with at least $2\delta$ edges, then $\mu_{f}(G) \ge \delta$.
\\
(ii) If $\delta> \frac{n-k}{2}$, then $\mu_{f}(G)>\frac{n-k}{2}$.
\end{Theorem}
\p If $\Delta(L) \le 2$, then
$f = \frac{1}{2} f_{E(L)}$ is a fractional matching, and so $\mu_{f}(G) \ge \frac{|E(L)|}{2} \ge s$.
This proves (i).
To show (ii), we observe that $\delta > 0$.
As $n \ge n-k$, by Theorem \ref{le1}(i) and Theorem \ref{th2}(i), we may assume that $\frac{n-k}{2}<\delta<\frac{n}{2}$.
If $\kappa(G) \ge 2$, then by Theorem \ref{le1}(ii), $G$ contains a cycle $C_s$ with $s = |E(C_s)| \ge 2\delta$.
By Theorem \ref{th2}(i),
$\mu_{f}(G) \geq\delta>\frac{n-k}{2}$. Hence we assume that $G$ contains a cut vertex $u$.
Let $G_{1}$ and $G_{2}$ are two of connected components of $G-u$. Let $v_{1}, v_{2} \in N_G(u)$
with $v_{1}\in V(G_{1})$ and $v_{2}\in V(G_{2})$. As $G$ is simple, $\delta(G-u) \ge \delta(G) - 1.$
In particular, $\delta(G_i) \ge \delta(G) - 1$, for $i\in \{1, 2\}$.
By Theorem \ref{le1}(iii), for $i \in \{1,2\}$, each $G_i$ contains a $(v_i, v_i')$-path $P_i$, for some $v_i' \in V(G_i - v_i)$,
with $|E(P_i)| \ge \delta-1,$ then $G$ contains a path $P = G[E(P_1) \cup v_{1}uv_{2}\cup E(P_2)]$ of length $2\delta$, and
by Theorem \ref{th2}(i), $\mu_{f}(G)\geq\delta>\frac{n-k}{2}$.
\q

The main tool in our paper is the following eigenvalue interlacing technique.
Given two non-increasing real sequences
$\theta_{1}\geq \theta_{2}\geq \cdots \geq \theta_{n}$ and
$\eta_{1}\geq \eta_{2}\geq \cdots \geq \eta_{m}$
with $n>m,$ the second sequence is said to {\bf interlace}
the first one if $\theta_{i}\geq \eta_{i}\geq\theta_{n-m+i}$
for $i=1, 2, \ldots, m.$
The interlacing is {\bf tight} if exists an integer $k\in[0, m]$
such that $\theta_{i}=\eta_{i}$ for $1\leq i\leq k$ and $\theta_{n-m+i}=\eta_{i}$ for
$k+1\leq i\leq m.$

Let $M$ be the following $n\times n$ matrix
\[
M=\left(\begin{array}{ccccccc}
M_{1,1}&M_{1,2}&\cdots &M_{1,m}\\
M_{2,1}&M_{2,2}&\cdots &M_{2,m}\\
\vdots& \vdots& \ddots& \vdots\\
M_{m,1}&M_{m,2}&\cdots &M_{m,m}\\
\end{array}\right),
\]
whose rows and columns are partitioned into subsets $X_{1}, X_{2},\ldots ,X_{m}$ of $\{1,2,\ldots, n\}$.
Let $M_{i,j}$ denote the submatrix (called a block) of $M$ by deleting the rows in $\{1,2,\ldots, n\} - X_i$ and deleting
the columns in $\{1,2,\ldots, n\} - X_j$.
The {\bf quotient matrix} $R$ of the matrix $M$ (with respect to the given partition)
is the $m\times m$ matrix whose entries are the
average row sums of the blocks $M_{i,j}$ of $M$. The partition is {\bf equitable}
if each block $M_{i,j}$ of $M$ has constant row (and column) sum.

\begin{Theorem}(Brouwer and Haemers \cite{BrHa09, Haem95})\label{le2}
Let $M$ be a real symmetric matrix. Then the eigenvalues of every
quotient matrix of $M$ interlace the ones of $M.$ Furthermore, if the
interlacing is tight, then the partition is equitable.
\end{Theorem}

\begin{Theorem}(Haynsworth \cite{Hayn59}, You et al. \cite{YLH})\label{le3}
Let $M$ be a partitioned matrix, and $R$ be its equitable quotient matrix.
Then the eigenvalues of the quotient matrix $R$ are eigenvalues of $M$.
Furthermore, if $M$ is a nonnegative matrix, then the spectral radius of the quotient matrix $R$ equals to the spectral
radius of $M$.
\end{Theorem}


Let $A=(a_{ij})$ and $B=(b_{ij})$ be two $n\times n$ matrices. Define  $A\leq B$
if $a_{ij}\leq b_{ij}$ for all $i$ and $j$; and $A<B$ if
$A\leq B$ and $A\neq B$.

\begin{Theorem}(Berman and Plemmons \cite{BePl79}, Horn and Johnson \cite{HoJo86})\label{le4}
Let $A=(a_{ij})$ and $B=(b_{ij})$ be two $n\times n$ matrices with the spectral radii $\lambda_{1}(A)$ and $\lambda_{1}(B)$. If $0\leq A\leq B$, then $\lambda_{1}(A)\leq \lambda_{1}(B)$. Furthermore, if $B$ is irreducible and $0\leq A<B$, then $\lambda_{1}(A)<\lambda_{1}(B)$.
\end{Theorem}

\section{Proof of Theorems \ref{th1} and \ref{th5}}

Following the notation in \cite{BoMu08}, for $S, T \subseteq V(G)$, define $E_G([S, T] = \{uv \in E(G): u \in S$ and
$v \in T\}$. We will use (\ref{f}) as the definition of the function $\varphi(a, n,\delta, k)$.
Throughout this section, we always assume that $n$ is a positive integer,
$a$ and $k$ are real numbers satisfying $a \ge 0$ and $0 < k < n$, and that
$G$ is a connected simple graph with $n = |V(G)|$,  $\delta = \delta(G) \le \frac{n-k}{2}$.

\subsection{Proof of Theorem \ref{th1}}

By contradiction, we assume that $\mu_{f}(G)\leq\frac{n-k}{2}$. By (\ref{BT}),
there exists a vertex subset $S\subseteq V(G)$ satisfying $i(G-S)-|S|\geq k$.
Let $T$ be the set of isolated vertices
in $G-S$, $|S|=s$ and $|T|=t$. Then $s+t\leq n$ and $t=i(G-S)\geq s+k$, and so
$s\leq\frac{n-k}{2}$. As $T$ is the set of isolated vertices
in $G-S$, we observe that $N_G(T) \subseteq S$, and so $s\geq \delta$.

Let $X = E_G[S,T]$, $H = G[X]$ be the bipartite subgraph $G$ induced by $X,$
and  $r=|E(H)|$. Then $r\geq t\delta$. Accordingly, the quotient matrix $R(aD(H)+A(H))$ of
$aD(H)+A(H)$ with respect to the partition $(S, T)$ becomes:
\[
R(aD(H)+A(H))=\left(\begin{array}{ccccccc}
\frac{ar}{s}&\frac{r}{s}\\
\frac{r}{t}&\frac{ar}{t}\\
\end{array}\right).
\]
As the characteristic polynomial of $R(aD(H)+A(H))$ is
$\displaystyle \lambda^{2}-a(\frac{r}{s}+\frac{r}{t})\lambda+(a^{2}-1)\frac{r^{2}}{st}=0$,
direct computation yields
\begin{eqnarray}\label{e2}
\nonumber
\lambda_{1}(R(aD(H)+A(H)))&=&\frac{1}{2}\{a(\frac{r}{s}+\frac{r}{t})+\sqrt{a^{2}(\frac{r}{s}+\frac{r}{t})^{2}-4(a^{2}-1)\frac{r^{2}}{st}}\}
\\&=&\frac{1}{2}\{a(\frac{r}{s}+\frac{r}{t})+\sqrt{(a^{2}-1)(\frac{r}{s}-\frac{r}{t})^{2}+(\frac{r}{s}+\frac{r}{t})^{2}}\}.
\end{eqnarray}
By Theorems \ref{le4} and \ref{le2}, to reach a contradiction to the assumption of Theorem \ref{th1},
it suffices to show that in each of the following cases, we always have
\begin{equation} \label{con}
\lambda_{1}(R(aD(H)+A(H))) \ge \varphi(a, n,\delta, k),
\end{equation}

\begin{Case} $a=0$.
\end{Case}

By (\ref{e2}), we have $\lambda_{1}(R(aD(H)+A(H)))=r\sqrt{\frac{1}{st}}$. Thus
\begin{eqnarray*}
\lambda_{1}(aD(G)+A(G))&\geq&\lambda_{1}(aD(H\cup(n-s-t)K_{1})+A(H\cup(n-s-t)K_{1}))\\
&=&\lambda_{1}(aD(H)+A(H))\geq\lambda_{1}(R(aD(H)+A(H)))\\
&=&r\sqrt{\frac{1}{st}}\geq \delta\sqrt{\frac{t}{s}}\geq \delta\sqrt{1+\frac{k}{s}}\\
&\geq&\delta\sqrt{1+\frac{2k}{n-k}}.
\end{eqnarray*}

\begin{Case} $0 < a \le 1$.
\end{Case}
This implies that $a^{2}-1 \le 0$, and so by (\ref{e2}),
\[
\lambda_{1}(R(aD(H)+A(H))) \geq\frac{1}{2}\{a(\frac{r}{s}+
\frac{r}{t})+\sqrt{(a^{2}-1)(\frac{r}{s}+\frac{r}{t})^{2}+(\frac{r}{s}+\frac{r}{t})^{2}}\}=a(\frac{r}{s}+\frac{r}{t}).
\]
It follows that
\begin{eqnarray}\label{e10}
\nonumber
\lambda_{1}(aD(G)+A(G))&\geq&\lambda_{1}(aD(H)+A(H))\geq\lambda_{1}(R(aD(H)+A(H)))
\\
&\geq& a(\frac{r}{s}+\frac{r}{t})\geq a\delta(1+\frac{t}{s})\geq a\delta(2+\frac{k}{s})
\ge  \frac{2a\delta n}{n-k}.
\end{eqnarray}

\begin{Case} $1 < a < +\infty$.
\end{Case}
This implies that $a^{2}-1>0$, and so by (\ref{e2}),
\[
\lambda_{1}(R(aD(H)+A(H))) >\frac{1}{2}\{a(\frac{r}{s}+\frac{r}{t})+
\sqrt{(a^{2}-1)(\frac{r}{s}-\frac{r}{t})^{2}+(\frac{r}{s}-\frac{r}{t})^{2}}\}=\frac{ar}{s}.
\]
It follows that
\begin{eqnarray*}
\lambda_{1}(aD(G)+A(G))&\geq&\lambda_{1}(aD(H)+A(H))\geq\lambda_{1}(R(aD(H)+A(H)))\\
&>&\frac{ar}{s}\geq a\delta\frac{t}{s}\geq a\delta(1+\frac{k}{s})
\ge \frac {a\delta (n+k)}{n-k}.
\end{eqnarray*}
As in all the possible values of $a$, (\ref{con}) is always obtained,
Theorem \ref{th1} is now justified.
\q

\subsection{Sharpness Discussion of Theorem \ref{th1}}

Following \cite{BoMu08}, for disjoint vertex sets $X$ and $Y$, $G[X,Y]$ denotes a bipartite graph with
bipartition $(X, Y)$.
We examine the following example to see that Theorem \ref{th1} is best possible when $a \in \{0,1\}$ in some sense.

\begin{Example}
Let  $X$ and $Y$ be disjoint vertex sets, and $k$, $n$ and $\delta$ be positive integers satisfying
$n > k > 0$ and $\delta \le \min\{|X|, |Y|\}$.
Define $\mathcal{B}(\delta, k)$ to be the family of bipartite graphs of the form $G[X, Y]$ with
$n = |X| + |Y|$, satisfying each of the following:
\\
(H1) for any $v \in X$, $d(v)$ is a constant,
\\
(H2) for any $v \in Y$, $d(v) = \delta$,
\\
(H3) $|Y|=|X|+k$.
\end{Example}
By definition, the complete bipartite graph $K_{\delta, \delta+k}\in \mathcal{B}(\delta, k)$.

\begin{Lemma}\label{le5}
For any $H^{*}\in\mathcal{B}(\delta, k)$, each of the following holds.
\\
(i) $\mu_{f}(H^{*})=\alpha'(H^{*})=|X|=\frac{n-k}{2}$.
\\
(ii) $a=0$ or $1$, $\lambda_{1}(aD(H^{*})+A(H^{*}))=\varphi(a, n,\delta, k).$
\\
(iii) $a\in(0,1)$ or $(1, +\infty)$, $\lambda_{1}(aD(H^{*})+A(H^{*}))>\varphi(a, n,\delta, k).$
\end{Lemma}
\p Let $d$ be the degree of each vertex in $X$.
By (H1) and (H2), $d|X| = |E(H^*)| = \delta|Y| = \delta(|X| + k)$, leading to
$d > \delta \ge 1$. For any subset $S \subseteq X$, $|N_{H^*}(S)| = d|S| > |S|$. It follows by Hall's Theorem
(Theorem 16.4 of \cite{BoMu08}) that $H^*$ has a matching $M$ of size $|X|$. Since $H^*$ is bipartite, by Theorem \ref{2mu}(ii),
we conclude that $\mu_{f}(H^{*})=\alpha'(H^{*})=|X|=\frac{n-k}{2}$, and so (i) follows.

We prove both (ii) and (iii) simultaneously.
The quotient matrix $R(aD(H^{*})+A(H^{*}))$ of $aD(H^{*})+A(H^{*})$ with respect to
the partitions $X$ and $Y$ of $V(H^*)$ can be written as
\begin{equation*}\begin{split}
R(aD(H^{*})+A(H^{*}))&=\left(\begin{array}{ccccccc}
\frac{a\delta|Y|}{|X|}&\frac{\delta|Y|}{|X|}\\[1mm]
\delta& a\delta\\
\end{array}\right).
\end{split}\end{equation*}
As the characteristic polynomial of $R(aD(H^{*})+A(H^{*}))$ is
$$\lambda^{2}-a\delta(1+\frac{|Y|}{|X|})\lambda+(a^{2}-1)\delta^{2}\frac{|Y|}{|X|}=0,$$
direct computation indicates that
\begin{eqnarray}\label{e6}
\nonumber
\lambda_{1}(R(aD(H^{*})+A(H^{*})))&=&\frac{1}{2}\{a\delta(\frac{|Y|}{|X|}+1)+\sqrt{a^{2}\delta^{2}(\frac{|Y|}{|X|}+1)^{2}-4(a^{2}-1)\delta^{2}\frac{|Y|}{|X|}}\}
\\&=&\frac{1}{2}\{a\delta(\frac{|Y|}{|X|}+1)+\sqrt{(a^{2}-1)\delta^{2}(\frac{|Y|}{|X|}-1)^{2}+\delta^{2}(\frac{|Y|}{|X|}+1)^{2}}\}.
\end{eqnarray}
As the partition is equitable and $aD(H^{*})+A(H^{*})$ is nonnegative, it follows by Theorem \ref{le3} that
$\lambda_{1}(aD(H^{*})+A(H^{*}))=\lambda_{1}(R(aD(H^{*})+A(H^{*})))$.
This, together with (\ref{e6}), leads to the following observations.
\\
(A) If $a=0$, then
$\lambda_{1}(A(H^{*}))=\delta\sqrt{\frac{|Y|}{|X|}}=\delta\sqrt{1+\frac{k}{|X|}}=\delta\sqrt{1+\frac{2k}{n-k}}=\varphi(0, n,\delta, k)$.
\\
(B) If $a=1$, then
$\lambda_{1}(D(H^{*})+A(H^{*}))=\delta(\frac{|Y|}{|X|}+1)=\frac{2n\delta}{n-k}=\varphi(1, n,\delta, k)$.
\\
(C) If $0 < a < 1$, then $a^{2}-1<0$, and so
\begin{eqnarray}
\nonumber
\lambda_{1}(aD(H^{*})+A(H^{*}))&>&\frac{1}{2}\{a\delta(\frac{|Y|}{|X|}+1)+\sqrt{(a^{2}-1)\delta^{2}(\frac{|Y|}{|X|}+1)^{2}+\delta^{2}(\frac{|Y|}{|X|}+1)^{2}}\}\\
\nonumber
&=&a\delta(\frac{|Y|}{|X|}+1)=\frac{2a\delta n}{n-k}=\varphi(a, n,\delta, k).
\end{eqnarray}
\noindent
(D) If $1 < a < +\infty$, then $a^{2}-1>0$, and so
\begin{eqnarray}
\nonumber
\lambda_{1}(aD(H^{*})+A(H^{*}))&>&\frac{1}{2}\{a\delta(\frac{|Y|}{|X|}+1)+
\sqrt{(a^{2}-1)\delta^{2}(\frac{|Y|}{|X|}-1)^{2}+\delta^{2}(\frac{|Y|}{|X|}-1)^{2}}\}
\\
\nonumber
&=&a\delta(\frac{|Y|}{|X|}+1)=\frac{a\delta(n+k)}{n-k}=\varphi(a, n,\delta, k).
\end{eqnarray}
These observations (A)-(D) complete the proof of the lemma.
\q

\begin{Remark}
By Lemma \ref{le5}, there exist graphs $H^{*}\in\mathcal{B}(\delta, k)$ with minimum degree $\delta$ such that
$\lambda_{1}(A(H^{*}))=\varphi(0, n,\delta, k)$, $\mu_{f}(H^{*})=\frac{n-k}{2}$ and $\lambda_{1}(D(H^{*})+A(H^{*}))=\varphi(1, n,\delta, k)$, $\mu_{f}(H^{*})=\frac{n-k}{2}$. Hence the upper bound in Theorem \ref{th1} is best possible when $a \in \{0, 1\}$.
\end{Remark}

\subsection{The Proof and Discussion of Theorem \ref{th5}}

By contradiction, assume that $\mu_{f}(G)\leq\frac{n-k}{2}$. By (\ref{BT}),
there exists a vertex subset $S\subseteq V(G)$ such that $i(G-S)-|S|\geq k$.
Let $T$ be the set of isolated vertices in $G-S$. Let $|S|=s$ and $|T|=t.$ Then $s+t\leq n$ and $t=i(G-S)\geq s+k.$
As $\cup_{v \in T}N_G(v) \subseteq S$, we observe that $s\geq\delta$, and so $t\geq s+k\geq\delta+k$.
As $K_{t}\cup(n-t)K_{1}$ is a spanning subgraph of $\overline{G}$, it follows by
Theorem \ref{le4} that
\begin{eqnarray}
\nonumber
\lambda_{1}(aD(\overline{G})+A(\overline{G}))&\geq& \lambda_{1}(aD(K_{t}\cup(n-t)K_{1})+A(K_{t}\cup(n-t)K_{1}))\\
\nonumber
&=&(a+1)(t-1)\geq(a+1)(\delta+k-1),
\end{eqnarray}
leading to a contradiction to the hypothesis. This completes the proof of  Theorem \ref{th5}.
\q

\begin{Example} \label{Ex-b}
By definition, the complete bipartite graph $K_{\delta, \delta+k}\in\mathcal{B}(\delta, k).$
Direct computation yields that $\lambda_{1}(aD(\overline{K_{\delta, \delta+k}})+A(\overline{K_{\delta, \delta+k}}))=(a+1)(\delta+k-1)$,
and $\mu_{f}(K_{\delta, \delta+k})=\mu(K_{\delta, \delta+k})=\delta=\frac{n-k}{2}$.
In this sense that the upper bound in Theorem \ref{th5} is best possible.
\end{Example}

\section{Corollaries of Theorems \ref{th1} and \ref{th5}}

Throughout this section, we assume that $a$,  $b$ and $k$ are real numbers with
$a\geq0$, $b>0$ and $k\in(0, n)$. Again, the function $\varphi(a, n,\delta, k)$ is defined
in (\ref{f}). We start by considering the matrix $J(G: a,b) = aD(G)+bA(G)$.
Theorems \ref{th1} and \ref{th5} have the following seemingly more general forms.

\begin{Corollary}\label{co1}
Let $G$ be a connected simple graph with $n = |V(G)|$ and $\delta = \delta(G) \leq\frac{n-k}{2}$. If
$\displaystyle \lambda_{1}(aD(G)+bA(G))<b\varphi(\frac{a}{b}, n, \delta, k)$, then $\mu_{f}(G)>\frac{n-k}{2}$.
\end{Corollary}
\p If $\mu_{f}(G)\leq\frac{n-k}{2}$, then
by Theorem \ref{th1}, $\lambda_{1}(aD(G)+A(G))\geq \varphi(a, n, \delta, k)$.
As $aD(G)+bA(G)=b(\frac{a}{b}D(G)+A(G))$, it follows that
$\displaystyle \lambda_{1}(aD(G)+bA(G))=b\lambda_{1}(\frac{a}{b}D(G)+A(G))\geq b\varphi(\frac{a}{b}, n, \delta, k)$, a contradiction.
\q

\begin{Corollary}\label{co2}
Let $G$ be a connected simple graph of order $n$ with $n = |V(G)|$ and $\delta = \delta(G) \leq\frac{n-k}{2}$. If
$\displaystyle\lambda_{1}(aD(\overline{G})+bA(\overline{G}))<(a+b)(\delta+k-1)$, then $\mu_{f}(G)>\frac{n-k}{2}$.
\end{Corollary}
\p
If  $\mu_{f}(G)\leq\frac{n-k}{2}$, then by Theorem \ref{th5}, $\lambda_{1}(aD(\overline{G})+A(\overline{G}))\geq(a+1)(\delta+k-1)$.
As $aD(\overline{G})+bA(\overline{G})=b(\frac{a}{b}D(\overline{G})+A(\overline{G}))$, it follows that $\displaystyle \lambda_{1}(aD(\overline{G})+bA(\overline{G}))=b\lambda_{1}(\frac{a}{b}D(\overline{G})+A(\overline{G}))\geq b(\frac{a}{b}+1)(\delta+k-1)=(a+b)(\delta+k-1)$, a contradiction.
\q

\vskip 0.3cm

Finally we observe that Corollary \ref{co3} follows by letting $a=b=1$ in Corollaries \ref{co1} and \ref{co2},
and Corollary \ref{co4} can be obtained from  Corollaries \ref{co1} and \ref{co2} by
setting $a=\alpha$ and $b=1-\alpha$.

\subsection{Fractional matching number and eigenvalues}

We present an application of Theorem \ref{th1} in this subsection by showing
a relationship between $\mu_{f}(G)$ and $\lambda_{1}(aD(G)+A(G))$ for a graph $G$.
Theorem \ref{th3}(i) was proved by O in \cite{SO16}, and is included here for the sake of completeness.

\begin{Theorem}\label{th3}
Let $G$ be a connected simple graph with $n = |V(G)|$ and $\delta = \delta(G) \leq\frac{n-k}{2}$. Each of the following holds.
\\
(i) (O \cite{SO16}) If $a=0$, then $\mu_{f}(G)\geq\frac{n\delta^{2}}{\lambda_{1}(G)^{2}+\delta^{2}}$,
where the equality holds if and only if $G\in \mathcal{B}(\delta, k)$ and
$k=\frac{n(\lambda_{1}(G)^{2}-\delta^{2})}{\lambda_{1}(G)^{2}+\delta^{2}}$ is an integer.
\\
(ii) If $a=1$, then $\mu_{f}(G)\geq\frac{\delta n}{\lambda_{1}(D(G)+A(G))}$,
where the equality holds if and only if $G\in \mathcal{B}(\delta, k)$ and $k=n-\frac{2\delta n}{\lambda_{1}(D(G)+A(G))}$ is an integer.
\\
(iii) If $0 < a < 1$, then $\mu_{f}(G)\geq\frac{a\delta n}{\lambda_{1}(aD(G)+A(G))}$. Furthermore,
the equality holds if $G\in \mathcal{B}(\delta, k)$ and $k=n-\frac{2a\delta n}{\lambda_{1}(aD(G)+A(G))}$ is an integer.
\\
(iv) If $1 < a < +\infty$, then $\mu_{f}(G)\geq\frac{a\delta n}{\lambda_{1}(aD(G)+A(G))+a\delta}$,
where the equality holds if $G\in \mathcal{B}(\delta, k)$ and
$k=\frac{n(\lambda_{1}(aD(G)+A(G))-a\delta)}{\lambda_{1}(aD(G)+A(G))+a\delta}$ is an integer.
\end{Theorem}
\p We only prove Theorem \ref{th3} (ii)-(iv).
Throughout the proof, we denote $\lambda_{1}(aD(G)+A(G))$ and $\mu_{f}(G)$ by $\lambda_{1}$ and $\mu_{f},$ respectively.
Theorem \ref{th1} says that if $\lambda_{1}<\varphi(a, n, \delta, k)$, then $\mu_{f}>\frac{n-k}{2}$.
Theorem \ref{th3} follows by applying Theorem \ref{th1} to each of the following possible values of the number $a$.

Assume that $0 < a \le 1$.
Define $k_1=n-\frac{2a\delta n}{\lambda_{1}}$.
As $\frac{1}{n-x}$ is an increasing function of $x$ on $(0, n)$, $\frac{2a\delta n}{n-k}$ decreases towards
$\lambda_{1}$ as $k$ decreases towards $k_1$, and
$\displaystyle \lim_{k \rightarrow k_1} \frac{2a\delta n}{n-k} = \lambda_1$.
By Theorem \ref{th1}, $\mu_{f}>\frac{n-k}{2}$ for each value of $k\in(k_1, n).$
It follows that $\displaystyle \mu_{f} \ge \lim_{k \rightarrow k_1} \frac{n-k}{2} = \frac{n-k_1}{2}=\frac{a\delta n}{\lambda_{1}}$.
Finally, if $G\in \mathcal{B}(\delta, k)$ and $k=n-\frac{2a\delta n}{\lambda_{1}(aD(G)+A(G))}$ is an integer, then direct computation yields
$\mu_{f}(G)=\alpha'(G)=\frac{n-k}{2}=\frac{a\delta n}{\lambda_{1}(aD(G)+A(G))}$. This proves Theorem \ref{th3}(iii) and
lower bound part of Theorem \ref{th3}(ii).

Now suppose that $a=1$ and that $\mu_{f}(G)=\frac{\delta n}{\lambda_{1}(D(G)+A(G))}$.
To simplify the notation and terms, we shall adopt the definition of the vertex sets $S$ and $T$, and the
definition of the graph $H$ in
Subsection 3.1, with $s = |S|$ and $t = |T|$.
As $\mu_{f}(G)=\frac{\delta n}{\lambda_{1}(D(G)+A(G))}$, every inequality in the arguments of the previous paragraph
must be an equality, and so $k=k_1$ and $\lambda_{1}=\frac{2\delta n}{n-k}$.
Furthermore, every inequality in (\ref{e10}) in Subsection 3.1
must also be an equality.
Hence in (\ref{e10}) we must have $\displaystyle a\delta(1+\frac{t}{s}) = a\delta(2+\frac{k}{s})$,
leading to $k = t-s$, and so $k$ is  an integer.
Likewise, in (\ref{e10}) we must also have $\displaystyle 2+\frac{k}{s} = \frac{2n}{n-k}$. This, together
with $k = t-s$, implies $n=s+t$.
It follows that $a=t\delta, t=k+s.$ By Theorem \ref{le2}, the partition is equitable, and thus $G \cong H \in \mathcal{B}(\delta, k)$.
This completes the proof of Theorem \ref{th3}(ii).

Finally we assume that $a\in(1, +\infty)$. Define $k_2=\frac{n(\lambda_{1}-a\delta)}{\lambda_{1}+a\delta}$.
As $\frac{n+x}{n-x}$ is an increasing function of $x$ on $(0, n)$,  $\frac{a\delta(n+k)}{n-k}$ decreases towards
$\lambda_{1}$ as $k$ decreases towards $k_2$ and
$\displaystyle \lim_{k \rightarrow k_2}  \frac{a\delta(n+k)}{n-k} = \lambda_1$.
By Theorem \ref{th1}, $\mu_{f}>\frac{n-k}{2}$ for each value of $k\in(k_2, n)$, and so
$\displaystyle \mu_{f} \ge \lim_{k \rightarrow k_2} \frac{n-k}{2} = \frac{n-k_2}{2}=\frac{a\delta n}{\lambda_{1}+a\delta}$.
If $G\in \mathcal{B}(\delta, k)$ with $k=\frac{n(\lambda_{1}(aD(G)+A(G))-a\delta)}{\lambda_{1}(aD(G)+A(G))+a\delta}$ being
an integer, then direct computation yields that
$\mu_{f}(G)=\alpha'(G)=\frac{n-k}{2}=\frac{a\delta n}{\lambda_{1}(aD(G)+A(G))+a\delta}$. This justifies Theorem \ref{th3}(iv),
and completes the proof of the theorem.
\q

\subsection{Fractional perfect matchings and eigenvalues}

The following Corollaries \ref{th4} and \ref{th7} are immediate consequences of Theorems \ref{th1} and \ref{th5},
respectively. Letting $k=1$ in Theorems \ref{th1} and \ref{th5}, respectively, we conclude that
the hypotheses of Corollaries \ref{th4} and \ref{th7} respectively imply that $\mu_{f}(G)>\frac{n-1}{2}$.
By Theorem \ref{2mu}, $2\mu_{f}(G)$ is an integer
and $\mu_{f}(G)\leq\frac{n}{2}$. This forces that $\mu_{f}(G)=\frac{n}{2}$, and so $G$ has a fractional
perfect matching.

\begin{Corollary}\label{th4}
Let $a$ be a real number with $a\geq0$, and $G$ be a simple connected graph of order $n$ with minimum degree $\delta\leq\frac{n-1}{2}.$ If
$$\lambda_{1}(aD(G)+A(G))<\varphi(a, n, \delta, 1),$$ then $G$ has a fractional perfect matching.
\end{Corollary}

\begin{Corollary}\label{th7}
Let $a$ be a real number with $a\geq0$, $G$ be a connected graph of order $n$ with minimum degree $\delta\leq\frac{n-1}{2},$ and $\overline{G}$ be the complement of $G$. If $$\lambda_{1}(aD(\overline{G})+A(\overline{G}))<(a+1)\delta,$$ then $G$ has a fractional perfect matching.
\end{Corollary}

As it is shown in Example \ref{Ex-b}, the upper bound in Corollary \ref{th7} is best possible.
This upper bound in Corollary \ref{th7} can be slightly improved by excluding some graphs.
Following \cite{BoMu08}, the {\bf join} of two graphs $G_1$ and $G_2$, denoted $G_1 \bigvee G_2$, is the graph formed from
the disjoint union of $G_1$ and $G_2$ by adding new edges joining
every vertex of $G_1$ to every vertex of $G_2$.

\begin{Theorem}
Let $a$ be a real number with $a\geq0$, $G$ be a connected graph with $n= |V(G)|$,
and $\delta = \delta(G) \leq\frac{n-1}{2}$. If
\begin{equation} \label{fpm}
\lambda_{1}(aD(\overline{G})+A(\overline{G}))<(a+1)(\delta+1),
\end{equation}
then $G$ has a fractional perfect matching unless $G\cong(\delta+1)K_{1}
\bigvee H_{\delta},$ where $H_{\delta}$ is a simple graph of order $\delta$.
\end{Theorem}
\p
By contradiction, assume that $G$ has no fractional perfect matchings.
By (\ref{T1}), there exists a vertex set $S\subseteq V(G)$ such that $i(G-S)-|S|>0$. Let $T$ be
the set of isolated vertices in $G-S$. Then $\cup_{v \in T} N_G(v) \subseteq S$, and so $|S|\geq \delta$.
This implies that $|T|\geq |S|+1\geq \delta+1$.

We claim that $V(G)=S\cup T$. By contradiction, assume that there exists a vertex $v\in V(G) -  S\cup T$.
Then we have $K_{\delta+2} \subseteq \overline{G}[T\cup \{v\}]\subseteq \overline{G}$.
By Theorem \ref{le4},
\begin{eqnarray}
\nonumber
\lambda_{1}(aD(\overline{G})+A(\overline{G}))&\geq&\lambda_{1}(aD(K_{\delta+2}\cup(n-\delta-2)K_{1})+A(K_{\delta+2}\cup(n-\delta-2)K_{1}))\\
&=&\lambda_{1}(aD(K_{\delta+2})+A(K_{\delta+2}))=(a+1)(\delta+1),
\nonumber
\end{eqnarray}
contrary to (\ref{fpm}).

Furthermore, we claim that $|T|=\delta+1$. Assume to the contrary that $|T|\geq\delta+2$. Then $\overline{G}$ contains
a clique of order $\delta+2$, and so $\lambda_{1}(aD(\overline{G})+A(\overline{G}))\geq(a+1)(\delta+1)$,
contrary to (\ref{fpm}) again. Hence $|T|=\delta+1$, and so $|S|=\delta$.
It follows that $G\cong (\delta+1)K_{1}\bigvee H_{\delta}$.
As there exists a vertex subset $V(H_{\delta})\subseteq V(G)$ such that
$i(G-V(H_{\delta}))>|V(H_{\delta})|$, by (\ref{T1}),
$(\delta+1)K_{1}\bigvee H_{\delta}$ has no fractional perfect matchings. This completes the proof.
\q

\vspace{5mm}
\noindent
{\bf Acknowledgement.}

The research of Ruifang Liu is supported by NSFC (No. 11971445), Natural Science Foundation of Henan Province (No. 20A110005).
The research of Hong-Jian Lai is supported by NSFC (Nos.~11771039 and 11771443).

\small {

}

\end{document}